\documentclass[12pt]{article}

\usepackage{amsfonts}

\newtheorem{thm}{Theorem}[section]
\newtheorem{lemma}[thm]{Lemma}
\newtheorem{cor}[thm]{Corollary}

\newtheorem{prop}[thm]{Proposition}
\newtheorem{conjecture}{Conjecture}

\newcommand{\proof
}{\par\medskip\noindent {\bf Proof.\ \ }}

\newcommand{\be}{\begin{equation}}
\newcommand{\ee}{\end{equation}}
\newcommand{\openbox}{\leavevmode
  \hbox to8pt{\hfil\vrule\vbox to6pt{\hrule width6pt\vfil\hrule}\vrule}}

\newcommand{\qed}{\hbox to5pt{ } \hfill \openbox\bigskip\medskip}

\newcommand{\cL}{\mbox{$\cal L$}}

\newcommand{\cF}{\mbox{$\cal F$}}
\newcommand{\cG}{\mbox{$\cal G$}}
\newcommand{\cH}{\mbox{$\cal H$}}

\newcommand{\cM}{\mbox{$\cal M$}}
\newcommand{\cK}{\mbox{$\cal K$}}

\title{A generalization of the Erd\H{o}s-Ko-Rado Theorem}
\author{G\'abor Heged\H{u}s
\\{\normalsize  \'Obuda University}
}
\begin{document}
\maketitle
\begin{abstract}
Our main result is a new upper bound for the size 
of $k$-uniform, $L$-intersecting families of sets, where $L$ contains only positive integers. 
We characterize extremal families in this setting. Our proof is based on the Ray-Chaudhuri--Wilson Theorem. 
As an application, we give a new proof for the Erd\H{o}s-Ko-Rado Theorem, improve Fisher's inequality in the uniform case and 
give an uniform version of the Frankl-F\"uredi conjecture .
\end{abstract}


\medskip
\section{Introduction}

First we introduce some notations.

Let  $[n]$ stand for the set $\{1,2,
\ldots, n\}$. We denote the family of all subsets of $[n]$  by $2^{[n]}$. 

For an integer $0\leq k\leq n$ we denote by
${[n] \choose k}$ the family of all  $k$ element subsets of $[n]$.

We call a family $\cF$ of subsets of $[n]$ {\em $k$-uniform}, if $|F|=k$ for each $F\in \cF$.
 
Bose  proved the following result in \cite{B}.
\begin{thm} \label{Bose}
Let $\lambda>0$ be a positive integer. Let $\cF=\{F_1,\ldots,F_m\}$ be a $k$-uniform family  of subsets of $[n]$ such that $|F_i\cap F_j|=\lambda$ for each $1\leq i,j\leq m, ~i \neq j$. Then $m\leq n$.
\end{thm}

Majumdar generalized this result and proved the following nonuniform version of Theorem \ref{Bose}. 
\begin{thm} \label{Maj}
Let $\lambda>0$ be a positive integer. Let $\cF=\{F_1,\ldots,F_m\}$ be a family of subsets of $[n]$ 
 such that $|F_i\cap F_j|=\lambda$ for each $1\leq i,j\leq m, ~i \neq j$. Then $m\leq n$.
\end{thm}

Frankl and F\"uredi conjectured in \cite{FF} and Ramanan proved in \cite{R} the following statement. 

\begin{thm} \label{Ram}
Let $\cF=\{F_1,\ldots,F_m\}$ be a family of subsets of $[n]$ 
 such that $1\leq |F_i\cap F_j|\leq s$ for each $1\leq i,j\leq m, ~i \neq j$. Then
$$
m\leq \sum_{i=0}^s {n-1 \choose i}.
$$
\end{thm}

Later Snevily conjectured the following statement in his doctoral dissertation (see \cite{S3}).  Finally he proved this result in \cite{S2}.

\begin{thm} \label{Snevily}
Let $\cF=\{F_1,\ldots,F_m\}$ be a family of subsets
of $[n]$. Let $L=\{\ell_1,\ldots ,\ell_s\}$ be a collection of $s$ positive integer. If $|F_i \cap F_j|\in L$ for each $1\leq i,j\leq m, ~i \neq j$, then 
$$
m\leq \sum_{i=0}^s {n-1 \choose i}.
$$
\end{thm}

A family $\cF$ is {\em $t$-intersecting}, if $|F\cap F'|\geq t$ whenever $F, F'\in \cF$. Specially, $\cF$ is an {\em intersecting} family, if $F\cap F'\ne \emptyset$ whenever $F, F'\in \cF$.

Erd\H{o}s, Ko and  Rado proved the following well-known result in \cite{EKR}:
\begin{thm} \label{EKZ}
Let $n,k,t$ be integers with $0<t<k< n$ . Suppose $\cF$ is a  $t$-intersecting,  $k$-uniform family of subsets of $[n]$. Then for $n>n_0(k,t)$,
$$
|\cF|\leq {n-t \choose k-t}.
$$
Further,  $|\cF|= {n-t \choose k-t}$ if and only if for some $T\in {[n]\choose t}$ we have 
$$
\cF=\{F\in {[n]\choose k}:~ T\subseteq F\}.
$$ 
\end{thm}
Let $L$ be a set of nonnegative integers. A family $\cF$ is {\em $L$-intersecting}, if $|E\cap F|\in L$ for every pair $E,F$ of distinct members of $\cF$. The following theorem gives a remarkable upper bound for the size of a $k$-uniform $L$-intersecting family. 

\begin{thm} \label{RW} (Ray-Chaudhuri--Wilson) Let $0<s\leq k\leq n$ be positive integers. Let $L$ be a set of $s$ nonnegative integers and $\cF=\{F_1,\ldots,F_m\}$  an $L$-intersecting, $k$-uniform family  of subsets of $[n]$. Then
$$
m\leq  {n \choose s}.
$$
\end{thm}

Erd\H{o}s, Deza  and Frankl improved Theorem \ref{RW} in \cite{EDF}. They used the theory of $\Delta$-systems in their proof.

\begin{thm} \label{EDF}
Let $0<s\leq k\leq n$ be positive integers. Let $L$ be a set of $s$ nonnegative integers and $\cF=\{F_1,\ldots,F_m\}$  an $L$-intersecting, $k$-uniform family  of subsets of $[n]$. Then for  $n>n_0(k,L)$
$$
m\leq \prod_{i=1}^s \frac{n-\ell_i}{k-\ell_i}.
$$ 
\end{thm}

Barg and Musin gave an improved version of Theorem \ref{RW} in \cite{BM}.
\begin{thm} \label{BM}
Let $L$ be a set of $s$ nonnegative integers and $\cF=\{F_1,\ldots,F_m\}$  an $L$-intersecting, $k$-uniform family  of subsets of $[n]$. Suppose that 
$$
\frac{s(k^2-(s-1))(2k-n/2)}{n-2(s-1)} \leq \sum_{i=1}^s \ell_i.
$$    
Then
$$
m\leq {n \choose s}-{n \choose s-1}\frac{n-2s+3}{n-s+2}.
$$
\end{thm}

First we prove a special case of our main result.
\begin{prop} \label{help}
Let $0<s\leq k\leq n$ be positive integers. Let $L=\{\ell_1 ,\ldots ,\ell_s\}$ be a set of $s$ positive integers such that $0<\ell_1 <\ldots <\ell_s$. Suppose that $n\geq {k^2 \choose \ell_1+1}s+\ell_1$. Let $\cF=\{F_1,\ldots,F_m\}$  be an $L$-intersecting, $k$-uniform family  of subsets of $[n]$. Suppose that $\bigcap\limits_{F\in \cF} \cF=\emptyset$. Then
$$
m\leq  {n-\ell_1 \choose s}.
$$
\end{prop}

We state now our main results.
\begin{thm} \label{main} 
Let $0<s\leq k\leq n$ be positive integers. 
 Let $L=\{\ell_1 ,\ldots ,\ell_s\}$ be a set of $s$ positive integers such that $0<\ell_1 <\ldots <\ell_s$. Suppose that $n\geq {k^2 \choose \ell_1+1}s+\ell_1$. Let $\cG=\{G_1,\ldots,G_m\}$  be an $L$-intersecting, $k$-uniform family  of subsets of $[n]$. Then
$$
m\leq  {n-\ell_1 \choose s}.
$$
\end{thm}

In the proof of Theorem \ref{main} we combine simple combinatorial arguments with the Ray-Chaudhuri--Wilson Theorem \ref{RW}. Our proof was inspired by the proof of Proposition 8.8 in \cite{J}.
\smallskip

In the following  we characterize the extremal families appearing in Theorem \ref{main}.
\begin{cor} \label{extrem}
Let $0<s\leq k\leq n$ be positive integers. 
Let $L=\{\ell_1 ,\ldots ,\ell_s\}$ be a set of $s$ positive integers such that $0<\ell_1 <\ldots <\ell_s$.  Suppose that $n> {k^2 \choose \ell_1+1}s+\ell_1$. Let $\cG=\{G_1,\ldots,G_m\}$  be an $L$-intersecting, $k$-uniform family  of subsets of $[n]$. Suppose that 
$$
|\cG|={n-\ell_1 \choose s}.
$$
Then there exists a $T\in {[n]\choose \ell_1}$ subset such that  $T\subseteq G$ for each $G\in \cG$. 
\end{cor}

We give here some immediate consequences of Theorem \ref{main}.
First we describe an uniform version of Theorem \ref{Ram}. 

\begin{cor} \label{main2}
Let $0<s<k\geq n$ be positive integers. Let $L=\{1 ,2,\ldots, s\}$. Suppose that $n>{k^2 \choose 2}s$. Let $\cF=\{F_1,\ldots,F_m\}$  be an $L$-intersecting, $k$-uniform family  of subsets of $[n]$. Then 
$$
m\leq  {n-1 \choose s}.
$$
Further if $n>{k^2 \choose 2}s+1$ and 
$$
|\cF|={n-1 \choose s},
$$
then $\bigcap_{F\in \cF} F\neq \emptyset$.
\end{cor}

The following result is the uniform version of Theorem \ref{Bose}.
\begin{cor} \label{Fisher}
Let $\lambda>0$ be a positive integer. Suppose that $n\geq{k^2 \choose \lambda+1}+\lambda$. Let $\cF=\{F_1,\ldots,F_m\}$ be a $k$-uniform family  of subsets of $[n]$ such that $|F_i\cap F_j|=\lambda$ for each $1\leq i,j\leq m, ~i \neq j$. Then
$$
m\leq n-\lambda.
$$
Further if $n>{k^2 \choose \lambda+1}+\lambda$ and 
$$
|\cF|= n-\lambda,
$$
then there exists a $T\in {[n]\choose \lambda}$ subset such that  $T\subseteq F$ for each $F\in \cF$.
\end{cor}
The following special case of Theorem \ref{EKZ}
follows immediately from Theorem \ref{main}.

\begin{cor} \label{EKZ2}
Let $n,k,t$ be integers with $0<t<k<n$. Suppose that $n\geq (k-t){k^2 \choose t+1}+t$. Let $\cF=\{F_1,\ldots,F_m\}$ be a  $t$-intersecting,  $k$-uniform family of subsets of $[n]$. Then 
$$
m\leq {n-t \choose k-t}.
$$
\end{cor}
\proof 
Let $L:=\{t, t+1, \ldots , k-1\}$ and apply Theorem \ref{main}.
\qed

Similarly Corollary \ref{extrem} implies the following result. 
\begin{cor} \label{EKZ3}
Let $0<k\leq n$ be integers such that $n> (k-t){k^2 \choose t+1}+t$. Let $\cF$ be a  $t$-intersecting,  $k$-uniform family of subsets of $[n]$. Suppose that 
$$
|\cF|={n-t \choose k-t}.
$$
Then there exists a $T\in{[n]\choose t}$ subset such that $T\subseteq F$ for each $F\in \cF$.
\end{cor}

\section{Proof}

The following Lemma is a well-known Helly-type result (see e.g. \cite{Bo}).
\begin{lemma} \label{Helly}
If each family of at most $k+1$ members of a $k$-uniform set system intersect, then all members intersect.
\end{lemma}
\qed

In our proof we use the following Lemma.

\begin{lemma} \label{soul} 
Let $\ell_1$ be a positive integer.
Let $\cH$ be a family  of subsets of $[n]$. Suppose that $\bigcap\limits_{H\in \cH} H =\emptyset$. 
Let $F\subseteq [n]$, $F\notin \cH$ be a subset such that $|F\cap H|\geq \ell_1$ for each $H\in \cH$. Let $Q:=\bigcup\limits_{H\in \cH} H$. Then
$$
|Q\cap F|\geq \ell_1+1.
$$
\end{lemma}
\proof
Since $|F\cap H|\geq \ell_1$ 
for each $H\in \cH$, thus $|Q\cap F|\geq \ell_1$. Indirectly, suppose that $|Q\cap F|= \ell_1$. Let $U:=Q\cap F$. Then
$$
U=Q\cap F=(\bigcup\limits_{H\in \cH} H)\cap F=\bigcup\limits_{H\in \cH} (H\cap F).
$$
Hence $H\cap F\subseteq U$ for each $H\in \cH$. Since $|U|=\ell_1$ and $|H\cap F|\geq \ell_1$ for each $H\in \cH$, thus $U=H\cap F$ for each $H\in \cH$. Hence $U\subseteq \bigcap\limits_{H\in \cH} H$, which is a contradiction with $\bigcap\limits_{H\in \cH} H =\emptyset$. \qed

\begin{lemma} \label{union}
 Let $\cH$ be a family  of subsets of $[n]$. Suppose that $t:=|\cH|\geq 2$ and $\cH$ is a $k$-uniform, intersecting family. Then
\begin{equation} \label{unio}
|\bigcup\limits_{H\in \cH} H|\leq k+(t-1)(k-1).
\end{equation}
\end{lemma}
\proof 
We use induction on  $t$. The inequality (\ref{unio}) is trivially true for $t=2$. 

Let $t\geq 3$. Suppose that the inequality  (\ref{unio}) is true for $t-1$. Let $\cH$ be an arbitrary $k$-uniform intersecting family such that $|\cH|=t$. Let $\cG\subseteq \cH$ be a fixed subset of $\cH$ such that $|\cG|=t-1$. Clearly $\cG$ is intersecting and $k$-uniform. It follows from the induction hypothesis that
$$
|\bigcup\limits_{G\in \cG} G|\leq k+(t-2)(k-1).
$$

Let $\{S\}=\cH \setminus \cG$. Then
$$
\bigcup\limits_{H\in \cH} H=(\bigcup\limits_{G\in \cG} G)\cup S, 
$$
thus 
$$
|\bigcup\limits_{H\in \cH} H|=|\bigcup\limits_{G\in \cG} G|+|S|-|(\bigcup\limits_{G\in \cG} G)\cap S|\leq k+(t-2)(k-1)+k-1=k+(t-1)(k-1).
$$
\qed

{\bf Proof of Proposition \ref{help}:}\\

Consider the special case when  $\bigcap\limits_{F\in \cF} F=\emptyset$.
By Lemma \ref{Helly} there exists a $\cG\subseteq \cF$ subset such that $\bigcap\limits_{G\in \cG} G=\emptyset$ and $|\cG|=k+1$. Let 
$$
M:=\bigcup\limits_{G\in \cG} G.
$$
It follows from Lemma \ref{union} that  $|M|\leq k+k(k-1)=k^2$. On the other hand it is easy to see that $|M\cap F|\geq \ell_1+1$ for each $F\in\cF$ by Lemma \ref{soul}. 

Let $T$ be a fixed subset of $M$ such that $|T|=\ell_1+1$.
Define
$$
\cF(T):=\{F\in \cF:~ T\subseteq M\cap F\}. 
$$
Let $L':=\{\ell_2, \ldots, \ell_s\}$. Clearly $|L'|=s-1$.
Then $\cF(T)$ is an $L'$-intersecting, $k$-uniform family, because $\cF$ is an $L$-intersecting family and $|M\cap F|\geq \ell_1+1$ for each $F\in\cF$. 
\begin{prop} \label{Funion}
$$
\cF=\bigcup\limits_{T\subseteq M, |T|=\ell_1+1} \cF(T).
$$
\end{prop}
\proof
Let  $\cM:=\bigcup\limits_{T\subseteq M, |T|=\ell_1+1} \cF(T)$. Clearly $\cM\subseteq \cF$. We prove that  $\cF\subseteq \cM$. \\
Let $F\in \cF$ be an arbitrary subset. Firstly, if $F\in \cG$, then $F\cap M=F$, because $M=\bigcup\limits_{G\in \cG} G$. Let $T$ be a fixed subset of $F$ such that $|T|=\ell_1+1$.  Then $F\in \cF(T)$.
Secondly, suppose that $F\notin \cG$. Then $|F\cap M|\geq \ell_1+1$ by Lemma  \ref{soul}. Let $T$ be a fixed subset of $F\cap M$ such that $|T|=\ell_1+1$.  Then  $F\in \cF(T)$ again. \qed

Let $T$ be a fixed, but arbitrary subset of $M$ such that $|T|=\ell_1+1$. Consider the set system
$$
\cG(T):=\{F\setminus T:~ F\in \cF(T)\}.
$$
Clearly $|\cG(T)|=|\cF(T)|$. Let $\overline{L}:=\{\ell_2-\ell_1-1, \ldots, \ell_s-\ell_1-1\}$. Here $|\overline{L}|=s-1$.
Since  $\cF(T)$ is an $L'$-intersecting, $k$-uniform family, thus $\cG(T)$ is an $\overline{L}$-intersecting, $(k-\ell_1-1)$-uniform family and $G\subseteq [n]\setminus T$ for each $G\in \cG(T)$. Hence it follows from Theorem \ref{RW} that
$$                   
|\cF(T)|=|\cG(T)|\leq {n-\ell_1-1 \choose s-1}.
$$
Finally Proposition \ref{Funion} implies that
$$
|\cF|\leq \sum_{T\subseteq M, |T|=\ell_1+1} |\cF(T)|\leq {k^2 \choose \ell_1+1} {n-\ell_1-1 \choose s-1}, 
$$
but
$$
{n-\ell_1-1 \choose s-1}=\frac{s}{n-\ell_1}{n-\ell_1 \choose s},
$$
hence
$$
|\cF|\leq {k^2 \choose \ell_1+1}\frac{s}{n-\ell_1}{n-\ell_1 \choose s} \leq {n-\ell_1 \choose s} 
$$
because $n\geq {k^2 \choose \ell_1+1}s+\ell_1$. \qed

{\bf Proof of Theorem \ref{main}:}\\

First we handle the case when  $|\bigcap\limits_{G\in \cG} G|\geq \ell_1$. Let $T$ be a fixed subset of $\bigcap\limits_{G\in \cG} G$ such that $|T|=\ell_1$. Consider the set system 
$$ 
\cK:=\{G\setminus T:~ G\in \cG\}.
$$ 
Obviously $|\cG|=|\cK|$.
Let $L':=\{0,\ell_2-\ell_1, \ldots ,\ell_s-\ell_1\}$. Then clearly $\cK$ is a ($k-\ell_1$)-uniform $L'$-intersecting set system   of subsets of $[n]\setminus T$. It follows immediately from Ray-Chaudhuri--Wilson Theorem \ref{RW} that 
$$
|\cG|=|\cK|\leq  {n-\ell_1 \choose s}.
$$

Now suppose that $|\bigcap\limits_{G\in \cG} G|=t$, where $0<t<\ell_1$. Let $T$ be a fixed subset of $\bigcap\limits_{G\in \cG} G$ such that $|T|=t$. Then consider the set system 
$$ 
\cF:=\{G\setminus T:~ G\in \cG\}.
$$ 
Clearly $|\cF|=|\cG|$.
Let $L':=\{\ell_1-t,\ell_2-t, \ldots ,\ell_s-t\}$. Then clearly $\cF$ is a ($k-t$)-uniform $L'$-intersecting set system   of subsets of $[n]\setminus T$. It follows from Proposition \ref{help} that 
$$
|\cG|=|\cF|\leq  {n-t-(\ell_1-t) \choose s}={n-\ell_1 \choose s}.
$$

Finally suppose that $\bigcap\limits_{G\in \cG} G=\emptyset$. Then Proposition \ref{help} gives us immediately that
$$
|\cG|\leq  {n-\ell_1 \choose s}.
$$
\qed

{\bf Proof of Corollary \ref{extrem}:}\\
It follows from the proof of Theorem \ref{main} that if $|\cF|= {n-\ell_1 \choose s}$ and $n>{k^2 \choose \ell_1+1}s+\ell_1$, then $|\bigcap\limits_{F\in \cF} F|\geq \ell_1$. 
Thus there exists a $T\in {[n]\choose \ell_1}$ such that $T\subseteq F$ for each $F\in \cF$.
\qed

\section{Remarks}

Let $q\geq 2$ stand for a fixed prime power. Let $PG(2,q)$ denote the finite projective plane over the Galois field $GF(q)$. Denote by $\cL$  the set of all lines of $PG(2,q)$. Let $k:=q+1$. Then $\cL$ can be considered as a $k$-uniform family of subsets of the base set $[k^2-k+1]$. Clearly $|\cL|=k$. 

This example motivates our next conjecture.
\begin{conjecture}
Let $0<s\leq k\leq n$ be positive integers. 
 Let $L=\{\ell_1 ,\ldots ,\ell_s\}$ be a set of $s$ positive integers such that $0<\ell_1 <\ldots <\ell_s$. Suppose that $n> k^2-k+1$. Let $\cF=\{F_1,\ldots,F_m\}$  be an $L$-intersecting, $k$-uniform family  of subsets of $[n]$. Then
$$
m\leq  {n-\ell_1 \choose s}.
$$
Further, if 
$$
|\cF|={n-\ell_1 \choose s},
$$
then there exists a $T\in {[n]\choose \ell_1}$ subset such that  $T\subseteq F$ for each $F\in \cF$.
\end{conjecture}

\end{document}